\documentclass{amsart}

\newtheorem{theorem}{Theorem}
\newtheorem{lemma}[theorem]{Lemma}

\newtheorem{corollary}[theorem]{Corollary}

\theoremstyle{definition}

\theoremstyle{remark}

\usepackage{amscd,amssymb}

\begin{document}

\title[Automorphisms fixing a variable of $K\langle x,y,z\rangle$]
{Automorphisms fixing a variable of $K\langle x,y,z\rangle$}

\author[Vesselin Drensky and Jie-Tai Yu]
{Vesselin Drensky and Jie-Tai Yu}
\address{Institute of Mathematics and Informatics,
Bulgarian Academy of Sciences,
          1113 Sofia, Bulgaria}
\email{drensky@math.bas.bg}
\address{Department of Mathematics, the University of Hong Kong,
          Hong Kong}
\email{yujt@hkusua.hku.hk}

\thanks
{The research of the first author was partially supported by Grant
MM-1106/2001 of the Bulgarian National Science Fund.}
\thanks{The research of the second author was partially
supported by RGC-CERG Grants 10203669 and 10204752.}

\subjclass{Primary 16S10. Secondary 16W20; 13B10; 13B25.}
\keywords{Automorphisms of free and polynomial algebras, tame automorphisms, wild
automorphisms}

\begin{abstract} We study automorphisms $\varphi$ of the free associative
algebra $K\langle x,y,z\rangle$ over a field $K$ such that
$\varphi(x),\varphi(y)$ are linear with respect to $x,y$ and
$\varphi(z)=z$. We prove that some of these automorphisms are wild
in the class of all automorphisms fixing $z$, including the well
known automorphism discovered by Anick, and show how to recognize
the wild ones. This class of automorphisms induces tame
automorphisms of the polynomial algebra $K[x,y,z]$. For $n>2$ the
automorphisms of $K\langle x_1,\ldots,x_n,z\rangle$ which fix $z$
and are linear in the $x_i$s are tame.
\end{abstract}

\maketitle

\section*{Introduction}

Let $K$ be a field of any characteristic and let
$X=\{x_1,\ldots,x_n\}$, $n\geq 2$, be a finite set. We denote by
$K[X]$ the polynomial algebra in the set of variables $X$ and by
$K\langle X\rangle$ the free associative algebra (or the algebra
of polynomials in the set $X$ of noncommuting variables). We write
the automorphisms of $K[X]$ and $K\langle X\rangle$ as $n$-tuples
of the images of the coordinates, i.e., $\varphi=(f_1,\ldots,f_n)$
means that $\varphi(x_j)=f_j(X)=f_j(x_1,\ldots,x_n)$,
$j=1,\ldots,n$. We distinguish two kinds of $K$-algebra
automorphisms of $K[X]$ and $K\langle X\rangle$. The first kind
are the affine automorphisms
\[
\left(\sum_{i=1}^n\alpha_{i1}x_i+\beta_1,\ldots,
\sum_{i=1}^n\alpha_{in}x_i+\beta_n\right),
\]
where $\alpha_{ij}$ and $\beta_j$ belong to $K$ and the $n\times n$ matrix
$(\alpha_{ij})$ is invertible. The second kind are the triangular automorphisms
\[
\left(\alpha_1x_1+f_1(x_2,\ldots,x_n),\ldots,
\alpha_{n-1}x_{n-1}+f_{n-1}(x_n),\alpha_nx_n+f_n\right),
\]
where $\alpha_j$ are invertible elements of $K$ and the
polynomials $f_j(x_{j+1},\ldots,x_n)$ do not depend on the first
$j$ variables. The affine and the triangular automorphisms
generate the group of the tame automorphisms. Instead of affine,
one can consider linear automorphisms only, assuming that the
polynomials $f_1,\ldots,f_n$ have no constant terms.

Many problems concerning automorphisms of free objects are stated in a
similar way for free groups, polynomial algebras, free associative and free Lie
algebras, for relatively free groups and algebras. Sometimes the
solutions are obtained with similar methods but very often they require
different techniques and the obtained results sound in different way,
see the recent book \cite{MSY} by Mikhalev, Shpilrain and Yu.
The results and the open problems on automorphisms of polynomial algebras
have served as the main motivation for the study of automorphisms of
free associative algebras.

It is a classical result of Jung \cite{J} (for $K=\mathbb C$) and
van der Kulk \cite{K} (for any $K$), that all automorphisms of
$K[x,y]$ are tame. Czerniakiewicz \cite{Cz} and Makar-Limanov
\cite{ML1, ML2} proved that all automorphisms of $K\langle
x,y\rangle$ are also tame. This implies that $\text{Aut }K\langle
x,y\rangle\cong \text{Aut }K[x,y]$.

It was an open problem
whether there exist nontame (or wild) automorphisms
of $K[X]$ for $n=\vert X\vert\geq 3$. Nagata \cite{N} constructed his famous
example, the automorphism of $K[x,y,z]$ defined by
\[
(x-2(y^2+xz)y-(y^2+xz)^2z,y+(y^2+xz)z,z).
\]
It fixes $z$ and, as Nagata showed, is wild considered as an
automorphism of $K[z][x,y]$, i.e., cannot be presented as a
product of tame automorphisms of $K[x,y,z]$ which fix $z$. It was
conjectured that the Nagata automorphism is wild also as an
element of $\text{Aut}(K[x,y,z])$ and it was one of the main open
problems in the theory of polynomial automorphisms for more than
30 years. Recently, Shestakov and Umirbaev \cite{SU1, SU2, SU3}
have developed a special technique based on noncommutative algebra
(Poisson algebras) and have established that the Nagata
automorphism is wild. They have also proved that every
automorphism of $K[x,y,z]$, which fixes $z$ and is wild as an
automorphism of $K[z][x,y]$, is wild as an automorphism of
$K[x,y,z]$. Umirbaev and Yu \cite {UY} have established even a
stronger version of this result: If $\varphi\in \text{\rm Aut
}K[z][x,y]$ is wild, then there exists no tame automorphism of
$K[x,y,z]$ which sends $x$ to $\varphi(x)$. In this way, if
$f(x,y,z)$ is a $K[z]$-wild coordinate in $K[z][x,y]$, then it is
immediately wild also in $K[x,y,z]$. (A coordinate means an
automorphism image of a variable.) In \cite{DY1} the authors of
the present paper commenced the systematic study of the wild
automorphisms and wild coordinates of $K[z][x,y]$, see also their
survey \cite{DY2}. In particular, they provided many new wild
automorphisms and wild coordinates of $K[z][x,y]$ which are
automatically wild automorphisms and wild coordinates of
$K[x,y,z]$. For $n>3$ the problem for existing of wild
automorphisms of $K[X]$ is still open.

Up till now no wild automorphisms of the free algebras with more
than two generators are known. There are some candidates to be
wild, see the book by Cohn \cite{C2}. One of them is the example
of Anick $(x+y(xy-yz),y,z+(zy-yz)y)\in \text{Aut }K\langle
x,y,z\rangle$, see \cite{C2}, p.~343. Although it fixes one
variable, its abelianization is a tame automorphism of $K[x,y,z]$
and we cannot apply the results of Shestakov, Umirbaev and Yu. Of
course, if we were able to lift the Nagata automorphism (or any
wild automorphism of $K[x,y,z]$) to any automorphism of $K\langle
x,y,z\rangle$, this would give an example of a wild automorphism
of $K\langle x,y,z\rangle$.

The present paper is motivated by the idea that the results on the
automorphisms of $K[z][x,y]$ give important information on the
automorphisms of $K[x,y,z]$. We study automorphisms of the free
algebra $K\langle x,y,z\rangle$ which fix the variable $z$. We
restrict our considerations to the automorphisms $\varphi$ such
that $\varphi(x),\varphi(y)$ are linear with respect to $x,y$. We
call such automorphisms linear $K[z]$-automorphisms. We prove that
some of these automorphisms are wild in the class of all
automorphisms fixing $z$, including the automorphism discovered by
Anick, and show how to recognize the wild ones. This class of
automorphisms induces tame automorphisms of the polynomial algebra
$K[x,y,z]$. For $n>2$ the automorphisms of $K\langle
x_1,\ldots,x_n,z\rangle$ which fix $z$ and are linear in the
$x_i$s are tame.

A forthcoming paper will be devoted to the description of the
structure of the group of all $K[z]$-automorphisms of $K\langle
x,y,z\rangle$.

\section{Preliminaries}

We fix a field $K$ of any characteristic, a set of variables
$X=\{x_1,\ldots,x_n\}$, $n\geq 2$, and one more variable $z$. The
main object of our paper is the free algebra $K\langle X,z\rangle$
in the set of free generators $X\cup\{z\}$. The algebra $K\langle
X,z\rangle$ is isomorphic to the free product $K[z]\ast_K K\langle
X\rangle$. We call an endomorphism $\varphi$ of $K\langle
X,z\rangle$ a $K[z]$-endomorphism if it fixes $z$ (and hence
$K[z]\subset K\langle X,z\rangle$) and write
$\varphi=(f_1,\ldots,f_n)$, where $f_j=\varphi(x_j)$,
$j=1,\ldots,n$. We denote the group of $K[z]$-automorphisms by
$\text{Aut}_{K[z]}K\langle X,z\rangle$.

Defining the tame $K[z]$-automorphisms of $K\langle X,z\rangle$,
the group of the triangular automorphisms consists obviously of
the automorphisms
\[
\left(\alpha_1x_1+f_1(x_2,\ldots,x_n,z),\ldots,
\alpha_{n-1}x_{n-1}+f_{n-1}(x_n,z),\alpha_nx_n+f_n(z)\right),
\]
$\alpha_j\in K^{\ast}$, $f_j\in K\langle X,z\rangle$, but we have
to decide which are the linear automorphisms.

We may call a $K[z]$-automorphism $(f_1,\ldots,f_n)$ linear if the
polynomials $f_j$ are linear with respect to $X$, with
coefficients depending on $z$. But for $n=2$ it is not known
whether a big class of automorphisms similar to this of Anick are
tame, although they all are linear with respect to $X$. So, we
prefer to introduce the group of elementary linear automorphisms
generated by the automorphisms
\[
(\alpha_1x_1,\ldots,\alpha_nx_n),\quad
(x_1,\ldots,x_{j-1},x_j+a(z)x_ib(z),x_{j+1},\ldots,x_n),
\quad i\not=j,
\]
$\alpha_j\in K^{\ast}$, $a,b\in K[z]$, and to generate the group
of the tame $K[z]$-automorphisms by the elementary linear
automorphisms and the triangular automorphisms. As we shall see,
the celebrated Suslin theorem \cite {Su} gives that for $n=\vert
X\vert \geq 3$ the two possible definitions are equivalent. Our
paper is concentrated around $K[z]$-endomorphisms of $K\langle
X,z\rangle$ which are linear with respect to $X$. We call such
endomorphisms linear $K[z]$-endomorphisms.

One of the main tools in the study of automorphisms is the
Jacobian matrix. If $\varphi=(f_1,\ldots,f_n)$ is an automorphism
of $K[X]$, then the Jacobian matrix $J(\varphi)=(\partial
f_j/\partial x_i)$ is invertible. The famous Jacobian conjecture
states that, if $\text{char }K=0$ and $J(\varphi)$ is invertible
for an endomorphism $\varphi$, then $\varphi$ is an automorphism,
see the book by van den Essen \cite{E} for the history and the
state-of-the-art of the problem. There are various analogues of
the Jacobian matrix in the case of (relatively) free groups and
(relatively) free associative and Lie algebras, see \cite{MSY}. In
the case of free associative algebras, the exact analogue of the
Jacobian matrix was introduced by Dicks and Lewin \cite{DL} who
proved the Jacobian conjecture for the free associative algebra
$K\langle x,y\rangle$ with two generators. The complete solution,
also into affirmative, was given by Schofield in his book
\cite{Sc}. Implicitly, in terms of endomorphisms, the Jacobian
matrix appeared in the paper \cite{Y} by Yagzhev who used it to
construct an algorithm which recognizes whether an endomorphism of
$K\langle X\rangle$ is an automorphism.

We recall the construction of the partial derivatives of Dicks and Lewin.
Let $y_0=z,y_1=x_1,\ldots,y_n=x_n$,
$Y=\{y_0,y_1,\ldots,y_n\}$ and $F=K\langle Y\rangle$. The algebra
${\mathcal M}(F)$ of the multiplications of $F$ is a subalgebra of the algebra of
$K$-linear operators acting on $F$ and is generated by the operators
$\lambda(u):F\to F$ and $\rho(u):F\to F$, $u\in F$,
of left and right multiplications defined, respectively, by
\[
w^{\lambda(u)}=uw,\quad w^{\rho(u)}=wu,\quad
w\in F.
\]
The operators $\lambda(u)$ and $\rho(v)$ commute
and ${\mathcal M}(F)\cong F^{\text{op}}\otimes_KF$,
where $F^{\text{op}}$ is the opposite algebra of $F$.
The isomorphism is given by
\[
\sum\lambda(u_p)\rho(v_p)\to
\sum u_p\otimes v_p
\]
and the opposite algebra appears because
$w^{\lambda(u_1)\lambda(u_2)}=u_2u_1w=w^{\lambda(u_2u_1)}$.
We define the partial derivatives on the monomials
$w=y_{j_1}\cdots y_{j_k}\in F$ by
\[
\frac{\partial w}{\partial y_i}=\sum_{p=1}^k
\left(y_{j_1}\cdots y_{j_{p-1}}\right)
\otimes\left(y_{j_{p+1}}\cdots y_{j_k}\right)
\delta_{pi}\in F^{\text{op}}\otimes_KF,
\]
where $\delta_{pi}=0,1$ is the Kronecker symbol, and then extend
them on $F$ by linearity. The Jacobian matrix of an endomorphism
$\varphi=(f_0,f_1,\ldots,f_n)$ of $F$ is defined by
\[
J(\varphi)=\left(\begin{matrix}
\frac{\partial f_0}{\partial y_0}&\frac{\partial f_1}{\partial y_0}&
\cdots&\frac{\partial f_n}{\partial y_0}\\
&&&\\
\frac{\partial f_0}{\partial y_1}&\frac{\partial f_1}{\partial y_1}&
\cdots&\frac{\partial f_n}{\partial y_1}\\
\vdots&\vdots&\ddots&\vdots\\
&&&\\
\frac{\partial f_0}{\partial y_n}&\frac{\partial f_1}{\partial y_n}&
\cdots&\frac{\partial f_n}{\partial y_n}\\
\end{matrix}\right).
\]
If $\varphi$ fixes $z=y_0$, the first column of this matrix consists of one 1
and $n$ zeros and the matrix is
invertible if and only if
the matrix
\[
J_{K[z]}(\varphi)=\left(\begin{matrix}
\frac{\partial f_1}{\partial y_1}&
\cdots&\frac{\partial f_n}{\partial y_1}\\
\vdots&\ddots&\vdots\\
\frac{\partial f_1}{\partial y_n}&
\cdots&\frac{\partial f_n}{\partial y_n}\\
\end{matrix}\right)
\]
is invertible and we call the latter matrix $J_{K[z]}(\varphi)$ the Jacobian
matrix of the $K[z]$-endomorphism $\varphi$.

It is well known (and easy to check) that the Jacobian matrices satisfy the chain rule
\[
J(\varphi\psi)=J(\varphi)\varphi(J(\psi)),
\]
where $\varphi(J(\psi))$ means that $\varphi$ acts on the entries
of $J(\psi)$. In particular, if $\varphi$ is an automorphism, then
$J_{\varphi}$ is invertible, i.e. belongs to $GL_{n+1}({\mathcal
M}(F))$. The theorem of Dicks-Lewin-Schofield gives that the
invertibility of $J(\varphi)$ implies that $\varphi$ is an
automorphism. It is obvious, that the same holds for the matrix
$J_{K[z]}(\varphi)$ when $\varphi$ is a $K[z]$-endomorphism of
$F=K\langle X,z\rangle$. If $\varphi=(f_1,\ldots,f_n)$ is a linear
$K[z]$-endomorphism of $K\langle X,z\rangle$, then
\[
f_j=\sum_{i=1}^n\sum_{p=1}^{k_{ij}}b_{ijp}(z)x_ic_{ijp}(z)
\]
and the entries $a_{ij}$ of $J_{K[z]}(\varphi)=(a_{ij})$ are of
the form
\[
a_{ij}=\sum_{p=1}^{k_{ij}}b_{ijp}(z)\otimes c_{ijp}(z)
\in K[z]^{\text{op}}\otimes_KK[z]\subset F^{\text{op}}\otimes_K F.
\]
Since $K[z]^{\text{op}}\otimes_KK[z]\cong K[z_1,z_2]$, identifying
$z\otimes 1$ with $z_1$ and $1\otimes z$ with $z_2$, we can
consider $J_{K[z]}(\varphi)$ as a matrix with entries from
$K[z_1,z_2]$.

\begin{lemma} \label{invertible Jacobian}
Let $\varphi$ be a linear $K[z]$-endomorphism of $F=K\langle
X,z\rangle$. Then $\varphi$ is an automorphism of $F$ if and only
if its Jacobian matrix $J_{K[z]}(\varphi)$ belongs to to the group
$GL_n(K[z_1,z_2])$ of the invertible matrices with entries from
$K[z_1,z_2]$. The group of the linear $K[z]$-automorphisms of $F$
is isomorphic to $GL_n(K[z_1,z_2])$.
\end{lemma}

\begin{proof}
The first part of the lemma follows from the fact that an $n\times
n$ matrix with entries from $K[z]^{\text{op}}\otimes_KK[z]\subset
F^{\text{op}}\otimes_KF$ is invertible over
$F^{\text{op}}\otimes_KF$ if and only if it is invertible over
$K[z]^{\text{op}}\otimes_KK[z]$. For the second part, if
$\varphi,\psi$ are linear $K[z]$-endomorphisms, then the matrices
$J_{K[z]}(\varphi), J_{K[z]}(\psi)$ depend on $z_1,z_2$ only and
the chain rule gives that
$J_{K[z]}(\varphi\psi)=J_{K[z]}(\varphi)J_{K[z]}(\psi)$, i.e., the
Jacobian matrix of the product of two linear $K[z]$-endomorphisms
is equal to the product, over $K[z_1,z_2]$, of the Jacobian
matrices of the factors.
\end{proof}

\section{The Main Results}

By the theorem of Suslin \cite{Su}, for $n\geq 3$, every matrix in
$GL_n(K[z_1,\ldots,z_p])$ can be presented as a product of a
diagonal matrix and elementary matrices, i.e., belongs to the
group $GE_n(K[z_1,\ldots,z_p])$. For $n=2$, this is not true. Cohn
\cite{C1} showed that the matrix
\[
\left(
\begin{matrix}
1+z_1z_2&z_2^2\\
z_1^2&1-z_1z_2\\
\end{matrix}
\right)\in SL_2(K[z_1,z_2])
\]
cannot be presented as a product of elementary $2\times 2$ matrices
with entries from $K[z_1,z_2]$.

\begin{theorem}\label{linear automorphisms}
{\rm (i)} A linear $K[z]$-automorphism of $K\langle x,y,z\rangle$
is tame if and only if its Jacobian matrix belongs to the group
$GE_2(K[z_1,z_2])$.
\par
{\rm (ii)} Every linear $K[z]$-automorphism of $K\langle x,y,z\rangle$
induces a tame automorphism of $K[x,y,z]$.
\par
{\rm (iii)} For $n>2$, any linear $K[z]$-automorphism of
$K\langle X,z\rangle$
is tame.
\end{theorem}

\begin{proof}
(i) The algebra $F=K\langle x,y,z\rangle$ has an augmentation
assuming that the variables $x,y$ are linear and $z$ is a
``constant'', i.e. of zero degree. The corresponding augmentation
ideal $\omega_F$ consists of all polynomials without terms
depending only on $z$. Every element $f$ of $F$ has the form
\[
f=f_0(z)+f_1(x,y,z)+f_2(x,y,z),
\]
where $f_0(z)\in K[z]$, $f_1(x,y,z)$ is linear with respect to
$x,y$ and $f_2(x,y,z)\in\omega^2_F$. It is easy to see (as in the
case of endomorphisms of $K[X]$ and $K\langle X\rangle$) that the
$K[z]$-endomorphism $\varphi=(f(x,y,z),g(x,y,z))$,
$f=f_0+f_1+f_2$, $g=g_0+g_1+g_2$, is an automorphism if and only
if the augmentation preserving endomorphism
$\varphi_0=(f_1+f_2,g_1+g_2)$ is an automorphism. Also, $\varphi$
is tame if and only if $\varphi_0$ is tame. Then we can decompose
$\varphi_0$ as a product of elementary augmentation preserving
automorphisms. So, we may restrict our considerations to
augmentation preserving $K[z]$-automorphisms only. If $\varphi_0$
is an automorphism, then the linear $K[z]$-endomorphism
$\varphi'=(f_1,g_1)$ is also an automorphism. If $\varphi_0$ is
tame, then $\varphi'$ is a product of elementary linear
$K[z]$-automorphisms. Hence, a linear $K[z]$-automorphism
$\varphi$ of $K\langle x,y,z\rangle$ is tame if and only if it is
a product of elementary linear $K[z]$-automorphisms. By Lemma
\ref{invertible Jacobian} this means that $J_{K[z]}(\varphi)$
belongs to $GE_2(K[z_1,z_2])$.

(ii) The linear $K[z]$-automorphism $\varphi$ of $K\langle
x,y,z\rangle$ induces a linear $K[z]$-auto\-morphism $\bar\varphi$
of $K[x,y,z]$. Its Jacobian matrix $J_{K[z]}(\bar\varphi)$ belongs
to $GL_2(K[z])$. Since $K[z]$ is a principal ideal domain, the
groups $GL_2(K[z])$ and $GE_2(K[z])$ coincide. This gives that
$J_{K[z]}(\bar\varphi)$ is a product of a diagonal matrix and
elementary matrices with entries from $K[z]$. Since the diagonal
and the elementary matrices correspond to elementary linear
automorphisms, we derive that $\bar\varphi$ is a tame
$K[z]$-automorphism.

(iii) The Jacobian matrix $J_{K[z]}(\varphi)$
of the linear $K[z]$-automorphism $\varphi$
of $K\langle X,z\rangle$ belongs to $GL_n(K[z_1,z_2])$.
The theorem of Suslin gives that $J_{K[z]}(\varphi)$ is a product of
a diagonal matrix and elementary matrices. Again, these matrices
correspond to elementary linear automorphisms and we obtain the proof.
\end{proof}

Recall that an automorphism $(f_1,\ldots,f_n)$ of $K\langle
X\rangle$ is called stably tame if the automorphism
$(f_1,\ldots,f_n,x_{n+1},\ldots,x_{n+m})$ of $K\langle
X,x_{n+1},\ldots,x_{n+m}\rangle$ is tame for some $m\geq 1$.
Theorem \ref{linear automorphisms} (iii) immediately gives:

\begin{corollary}\label{linear automorphisms are stably tame}
The linear $K[z]$-automorphisms of
$K\langle x,y,z\rangle$ are stably tame.
\end{corollary}

There is an algorithm which decides whether a matrix
in $GL_2(K[z_1,\ldots,z_p])$ belongs to $GE_2(K[z_1,\ldots,z_p])$.
It was suggested by Tolhuizen, Hollmann and Kalker \cite{THK}
for the partial ordering by degree and then generalized  by Park \cite{P}
for any monomial ordering on $K[z_1,\ldots,z_p]$. One applies Gaussian
elimination process on the matrix based on the Euclidean division algorithm
for $K[z_1,\ldots,z_p]$. The matrix belongs to $GE_2(K[z_1,\ldots,z_p])$
if and only if this procedure brings it to an elemetary or diagonal matrix.
The result of Park was already used by Shpilrain and Yu \cite{SY} to
give an algorithm which recognizes whether a polynomial in $K[x,y]$ is a coordinate,
and by the authors in \cite{DY1} to decide whether
a polynomial in $K[z][x,y]$ is a coordinate and a tame coordinate.

Consider the automorphism $(x+y(xy-yz),y,z+(zy-yz)y)$ of $K\langle
x,y,z\rangle$ constructed by Anick. Exchanging the places of $y$
and $z$, we obtain the automorphism
\[
\varphi=(x+z(xz-zy),y+(xz-zy)z)\in \text{\rm Aut}_{K[z]}K\langle x,y,z\rangle.
\]
Its abelianization $\bar\varphi=(x+z^2(x-y),y+z^2(x-y))$ is a tame
automorphism of $K[z][x,y]$: Apply Theorem \ref{linear
automorphisms} (ii) or change the coordinates $(x,y)$ of
$K[z][x,y]$ to $(u,y)=(x-y,y)$, then
\[
\bar\varphi=(x-y,y+z^2(x-y))=(u,y+z^2u).
\]
Clearly, $\varphi$ is a linear $K[z]$-automorphism. Its Jacobian matrix is
\[
J_{K[z]}(\varphi)=\left(
\begin{matrix}
1+z_1z_2&z^2\\
z_1^2&1-z_1z_2\\
\end{matrix}
\right)
\]
and is the matrix constructed by Cohn.
It cannot be presented as a product of elementary $2\times 2$ matrices
with entries from $K[z_1,z_2]$.
(Direct arguments: Applying the Gaussian
elimination process, we cannot reduce the entries of $J_{K[z]}(\varphi)$
using only the Euclidean division algorithm because the
leading terms of the entries of the columns are not divisible by each other.)
Hence, Theorem \ref{linear automorphisms} (i)
gives that this automorphism is wild considered as a $K[z]$-automorphism.
On the other hand, by Corollary \ref{linear automorphisms are stably tame},
it is stably tame. An explicit decomposition of $J_{K[z]}(\varphi)$
can be obtained from the proof of the lemma of Mennicke, see Lemma 2.3 of \cite{PW}.
The sequence of elementary operations in \cite{PW}, p. 281, formula (2.1),
gives the decomposition
\[
\left(
\begin{matrix}
1+z_1z_2&z_2^2&0\\
z_1^2&1-z_1z_2&0\\
0&0&1\\
\end{matrix}
\right)
\]
\[
=E_{13}(-z_2)E_{23}(-z_1)E_{31}(z_1)
E_{32}(-z_2)E_{13}(z_2)E_{23}(z_1)E_{31}(-z_1)E_{32}(z_2),
\]
where the matrix $E_{ij}(\alpha z_1^az_2^b)$ corresponds to the
$K[z]$-automorphism of $K\langle x,y,t,z\rangle$ defined by
$x_j\to x_j+\alpha z^ax_iz^b$, $x_k\to x_k$, when $k\not=j$, and
$x_1=x$, $x_2=y$, $x_3=t$.

\section*{Acknowledgements}

This project was carried out when the first author visited the
Department of Mathematics of the University of Hong Kong. He is
very grateful for the hospitality and the creative atmosphere
during his stay in Hong Kong.

\end{document}